\documentclass{amsart}
\usepackage{amsthm,amssymb}

\newtheorem{theorem}{Theorem}

\newtheorem{lemma}[theorem]{Lemma}

\newtheorem{remark}{Remark}

\newtheorem{conjecture}{Conjecture}

\begin{document}
\title[Polynomials ]{An upper bound on Jacobi polynomials}
\author[I. Krasikov]{Ilia Krasikov}
\address{ Department of Mathematical Sciences, Brunel University, Uxbridge
UB8 3PH United Kingdom} \email{mastiik@brunel.ac.uk}

\begin{abstract}
Let ${\bf P}_k^{( \alpha , \beta )} (x)$ be an orthonormal Jacobi
polynomial of degree $k.$ We will establish the following
inequality
\begin{equation*}
\max_{x \in [\delta_{-1},\delta_1]}\sqrt{(x-
\delta_{-1})(\delta_1-x)} \, (1-x)^{\alpha}(1+x)^{\beta}
\left({\bf P}_{k}^{( \alpha , \beta)} (x)\right)^2 < \frac{3
\sqrt{5}}{5} \, ,
\end{equation*}
where $\delta_{-1}<\delta_1$ are appropriate approximations to the
extreme zeros of ${\bf P}_k^{( \alpha , \beta )} (x) .$  As a
corollary we confirm, even in a stronger form,
 T. Erd\'{e}lyi, A.P. Magnus and P.
Nevai conjecture [Erd\'{e}lyi et al., Generalized Jacobi weights,
Christoffel functions, and Jacobi polynomials, SIAM J. Math. Anal.
25 (1994), 602-614], by proving that
\begin{equation*}
\max_{x \in
[-1,1]}(1-x)^{\alpha+\frac{1}{2}}(1+x)^{\beta+\frac{1}{2}}\left({\bf
P}_k^{( \alpha , \beta )} (x) \right)^2 < 3 \alpha^{1/3} \left(1+
\frac{\alpha}{k} \right)^{1/6} ,
\end{equation*}
 in the region $k
\ge 6, \; \; \alpha , \beta \ge \frac{1+ \sqrt{2}}{4}.$

\hspace{1pt}

\end{abstract}

\maketitle


\section{Introduction}
In this paper we will use bold letters for orthonormal polynomials
versus  regular characters for orthogonal polynomials in the
standard normalization \cite{szego}.

One of the most surprising and profound features of many families
of orthogonal polynomials is their equioscillatory behaviour. This
phenomenon has been discovered by G.\,Szeg\"{o} who proved that
for a vast class of weights the function $\sqrt{1-x^2} \,
\mathcal{W}(x) \, {\bf p}_i^2(x), $ asymptotically, for $i
\rightarrow \infty , $ equioscillates between $\pm \,
\frac{2}{\pi} \, .$ Here $\{ {\bf p}_i(x) \}$ is a family of
orthonormal polynomials of degree $i$ orthogonal with respect to
the weight function $\mathcal{W}(x)$ on $[-1,1],$ \cite{szego}.

A powerful theory developed for exponential weights
$\mathcal{W}=e^{-Q(x)}$ by A.L. Levin and D.S. Lubinsky
\cite{levlub} shows that under some smoothness assumptions on $Q,$
\begin{equation}
\label{eqlelu}
 \max_{I} \left| \sqrt{|(x- a_{-1})(a_1-x)|}
\, \mathcal{W}(x) \, {\bf p}_i^2(x) \right|< C,
\end{equation}
where the constant $C$ is independent on $i$ and $a_{\pm 1}=a_{\pm
1}(i)$ are Mhaskar-Rahmanov-Saff numbers for $Q.$  Recently
results of this type has been obtained for the Laguerre-type
exponential weights $x^{2 \rho}e^{-2 Q(x)}$ \cite{kasuga,levin4}.

It seems that in many cases (\ref{eqlelu}) is sharp, what means
that under an appropriate scaling, the envelope of the function
$\sqrt{|(x- \delta_{-1})(\delta_1-x)|}\; \mathcal{W}(x) {\bf
p}_i^2(x) ,$ where $\delta_{\pm 1}=\delta_{\pm 1}(i)$ are certain
approximations to the extreme zeros of $p_i$, is almost
independent on $i$ and has a plateau in the oscillatory region
with rapidly decaying slopes outside.

Nevertheless, today we don't possess even a vague picture which
could help us to expect such a behaviour. Already the case of
classical orthogonal polynomials is not properly covered by the
general theory. For example, two-side analogue of (\ref{eqlelu})
with explicit constants and valid independently not only on the
degree but also on all the parameters involved, are known for
Hermite \cite{krasherm} and Laguerre \cite{krasl2} polynomials.
Similar and, in a sense, best possible upper bound was also given
for the Bessel function \cite{landau}. However, it is not known
whether this is a peculiar property of the hypergeometric function
or a manifestation of a more general phenomenon.

 Surprisingly
enough, the first non-asymptotic inequality of this type was
obtained in the seemingly most complicated Jacobi case. Let
$$M_k^{\alpha , \beta }(x) = (1-x)^{\alpha+\frac{1}{2}}(1+x)^{\beta+\frac{1}{2}}\left({\bf
P}_k^{( \alpha , \beta )} (x) \right)^2 ,$$
 T.\,Erd\'{e}lyi, A.P.\,Magnus and P.\,Nevai,
proved that for  $k \ge 0, \; \; \alpha, \beta \ge - \frac{1}{2}
\, , $ \cite{erdel},
\begin{equation}
\label{eqernt}
 \mathcal{M}_k^{ \alpha , \beta}  = \max_{x \in
[-1,1]}M_k^{\alpha , \beta }(x) \le
 \frac{2e \left(2+ \sqrt{\alpha^2+\beta^2} \right)}{\pi} \, .
\end{equation}

They also observed that $\mathcal{M}_0^{ \alpha ,
\beta}=\frac{\sqrt{\alpha+ \beta}}{2 \pi} (1+o(1))$ for large
$\alpha, \beta$ and conjectured that the real order of
$\mathcal{M}_k^{ \alpha , \beta}$ is  $O \left( \alpha^{1/2}
\right).$

The inequality (\ref{eqernt}) was improved to
$$\mathcal{M}_k^{ \alpha , \beta} =O \left( \alpha^{2/3} \left(1+ \frac{\alpha}{k}
\right)^{1/3} \right), \; \; \alpha \ge \beta,$$ in
\cite{krasjac}, where a stronger conjecture was suggested,
 \begin{equation}
 \label{dokconj}
 \mathcal{M}_k^{ \alpha
, \beta}=O \left( \alpha^{1/3} \left(
1+\frac{\alpha}{k}\right)^{1/6} \right).
\end{equation}

 Under the classical restrictions $- \,\frac{1}{2} \le \alpha ,
\beta \le \frac{1}{2} \, ,$ much sharper inequalities implying
$$
\mathcal{M}_k^{ \alpha , \beta} \le \frac{2}{\pi} +O( \frac{1}{k}
),
$$
 are known \cite{chow} and \cite{lorch}.

Recently the author proved \cite{krasultsp} that in the
ultraspherical case for $ \alpha \ge \frac{1}{2} \, ,$

\begin{equation*}
\label{gleq2}  \max_{|x| \le \delta} \sqrt{\delta^2-x^2} \,
(1-x)^{\alpha }(1+x)^{\beta}{\bf P}_k^{ \alpha, \alpha }( x) <
\frac{2\nu}{\pi} \left(1+ \frac{1}{2(k+ \alpha)^2 }\right),
\end{equation*}
where
$$\delta= \sqrt{1- \frac{4\alpha^2-1}{(2k+2
\alpha+1)^2-4}} \, , \; \; \;  \nu = \left\{
\begin{array}{cc}
1 \, , & k \ge 4, \; \; even,\\
& \\
228, & k \ge 3, \; \; odd.
\end{array}
\right.
$$
This implies that for $\alpha \ge \frac{1+ \sqrt{2}}{4} \, , \; \;
k \ge 6,$
\begin{equation*}
\label{gleq1}
 \mathcal{M}_k^{ \alpha, \alpha } <21 \alpha ^{1/3}
\left(1+ \frac{\alpha }{k} \right)^{1/6}.
\end{equation*}

 The aim of this paper is to establish similar results in the
general case. In particular, we will confirm the above conjectures
under some mild restrictions on $ \alpha, \beta $ and $k.$ Namely
we prove the following:
\begin{theorem}
\label{thdelta} Let
\begin{equation}
\label{defz}
 Z_k^{\alpha, \beta}(x)=
\sqrt{\sqrt{(x-\delta_{-1})(\delta_1-x)} \; (1-x)^{ \alpha}
(1+x)^{ \beta}} \; {\bf P}_k^{(\alpha , \beta ) }(x),
\end{equation}
where
$$
\delta_{\pm 1}=\frac{\beta^2-\alpha^2 \pm \sqrt{(2k+1)(2k+2
\alpha+1)(2k+2 \beta+1)(2k+2 \alpha+2 \beta+1) }}{(2k+\alpha+
\beta+1)^2} \, .
$$
Then
\begin{equation}
\label{eqthdelta}
 \max_{x \in
[\delta_{-1},\delta_1]}|Z_k^{\alpha, \beta}(x)|
 <
5^{-1/4}\, \sqrt{3} ,
\end{equation}
 provided $k \ge 1, \; \; \alpha \ge \beta \ge 0.$
\end{theorem}
\begin{theorem}
\label{themn}
\begin{equation}
\label{eqthemn}
 \mathcal{M}_k^{ \alpha , \beta} < 3
\alpha^{1/3} \left(1+ \frac{\alpha}{k} \right)^{1/6},
\end{equation}
 provided
 $k \ge 6, \; \;  \alpha \ge \beta \ge
\frac{1+ \sqrt{2}}{4}.$
\end{theorem}
Notice that the interval $[\delta_{-1},\delta_1]\subset [-1,1]$
defined in Theorem \ref{thdelta} is large enough and contains, for
example, all the zeros of $P_k^{( \alpha , \beta )} (x).$
Probably, similar results hold for $ \alpha , \beta \ge - \,
\frac{1}{2} \, .$ Notice also that the assumption $ \alpha \ge
\beta$ does not impose any farther restrictions as $P_k^{ \alpha ,
\beta } (x)=(-1)^k P_k^{ \beta , \alpha } (-x).$

 We believe that
the bounds of Theorems \ref{thdelta} and \ref{themn} are sharp up
to a multiplicative factor. This means that there is a constant
$c$ such that, say, under the assumption of Theorem \ref{thdelta},
the inequality $|Z_k^{\alpha, \beta}(x)|
>c,$
holds in a point between any two consecutive zeros of
$P_k^{(\alpha , \beta ) }(x).$
 Furthermore, we suggest the following stronger conjecture
\begin{conjecture}
 \label{conjik1}
$$\max_{x \in [\delta_{-1},\delta_1]}
\sqrt{|(x-\delta_{-1})(\delta_1-x)|} \, (1-x)^{ \alpha} (1+x)^{
\beta}\left({\bf P}_k^{(\alpha , \beta ) }(x)\right)^2
=\frac{2}{\pi} \left(1+o(1) \right).$$
  provided
$ \alpha \ge \beta  \ge - \, \frac{1}{2} \,,$ and $(k+ \alpha )
\rightarrow \infty .$
\end{conjecture}
There is a good reason to believe that $P_k^{(\alpha , \beta )
}(x)$ lives on $[\delta_{-1}, \delta_1],$ namely
\begin{conjecture}
 \label{conjik2}
$$\int_{\delta_{-1}}^{\delta_1}
(1-x)^{\alpha}(1+x)^{\beta}\left({\bf P}_k^{( \alpha , \beta )}
(x) \right)^2 =1-o(1),$$
 for $ \alpha \ge \beta  \ge - \,
\frac{1}{2} \,,$ and $(k+ \alpha ) \rightarrow \infty .$
 \end{conjecture}

 Let us outline the proof of Theorems \ref{thdelta} and \ref{themn}.
 It is not difficult to find a pointwise upper bound on
$P_k^{\alpha, \beta}(x)$ in the bulk of the oscillatory region and
one such a bound has been already given in \cite{erdel}. In
section \ref{secosc} we establish a similar inequality which will
be enough for our purposes. Unfortunately, it seems that estimates
of this type diverge or become very poor in the transition region
around the extreme zeros. The damping factor
$\sqrt{(x-\delta_{-1})(\delta_1-x)} \,$ in (\ref{defz}) is needed
to move the global extremum into the oscillatory region. To prove
this we use so-called Sonin's function $S(f(x);x)$ which may be
viewed as an envelop of $f(x).$ Then it is a matter of simple
algebra to show that for the above choice of $\delta_{\pm 1}$ all
local extrema of $\left( Z_k^{\alpha , \beta }(x) \right)^2$ lay
on a curve with the unique maximum inside a proper subinterval of
$[ \delta_{-1}, \delta_1].$ This will be accomplished in section
\ref{secunim}. Theorems \ref{thdelta} will be proven in the last
section.  Passage to Theorem \ref{themn} is quite straightforward
since all the maxima of $M_k^{ \alpha , \beta } (x) $ belong to $[
\delta_{-1} +\epsilon, \delta_1-\epsilon],$ where $ \epsilon =O
\left( (k+ \alpha)^{-2/3} \right).$ This follows from the results
of \cite{krasjac}.

\section{Unimodality of Sonin's function}
\label{secunim}
 Jacobi polynomials $P_k^{( \alpha , \beta )} (x)$
are orthogonal polynomials on [-1,1] with the weight function
$\mathcal{W}(x)=(1-x)^{\alpha} (1+x)^{\beta},$ $\alpha, \beta
>-1,$ satisfying the following differential equation
\begin{equation}
\label{difjp} (1-x^2)y''+( \beta - \alpha -(\alpha + \beta
+2)x)y'+k (k+\alpha+ \beta+1) y, \; \; \; y=P_k^{( \alpha , \beta
)} (x).
\end{equation}
In the standard normalization their norm ${\bf h}_k$ is given by
\begin{equation}
\label{normjc} {\bf h}_k^2= \frac{2^{\alpha+ \beta+1} \Gamma (k+
\alpha+1) \Gamma (k+ \beta+1)}{(2k+ \alpha+ \beta+1) k! \, \Gamma
(k+ \alpha+ \beta+1)} \, .
\end{equation}
To avoid unnecessary discussion of some degenerate cases in the
sequel we will assume that
\begin{equation}
\label{assum1}
  \alpha > \beta >0.
\end{equation}
Theorems
\ref{thdelta} and \ref{themn} for $\alpha=\beta ,$ as well as for
$\beta=0,$ follow as obvious limiting cases.

To simplify otherwise messy expressions we will use the following
notation:
\begin{equation}
\label{note}
  \eta = \alpha - \beta , \; \; \; \sigma = \alpha+
\beta , \; \; \; r =2k+\alpha+ \beta +1,
\end{equation}
$$
 q=\eta/r \, , \; \; \;  s= \sigma/r  \, ;
$$
and also their trigonometric counterparts:
\begin{equation}
\label{notetr} q= \sin { \omega } , \; \; \;  s = \sin {\tau } .
\end{equation}
Thus, (\ref{assum1}) yields
\begin{equation}
\label{assum2}  0 < q <s <1, \; \; \; 0 < \omega <\tau <
\frac{\pi}{2} \,  .
\end{equation}

Using this notation we can rewrite $\delta_{\pm 1}$ defined in
Theorem \ref{thdelta} as follows:
 $\delta_j= j \cos{(\tau+j \,
\omega )} , \; \; \; j= \pm 1.$

We also introduce the function
$$
d(x)=(x- \delta_{-1} )(\delta_1-x)=1-q^2-s^2- 2 q  s x-x^2 .
$$
 We will see that in some respects the new
variables $q, s$ and $r$ are more natural than $ \alpha , \beta$
and $k.$

We start with the following simple lemma established in
\cite{krasultsp}. The proof is straightforward and is given here
for self-completeness.

 Given a real function $f(x),$ Sonin's function $S=S(f;x)$ is
$S=f^2+ \psi (x){f'^2},$ where $\psi(x) >0$ on an interval
$\mathcal{I}$ containing all local extrema of $f.$
 Thus, they lay on $S,$ and if $S$ is unimodal we can locate the global one.

\begin{lemma}
\label{mainl} Suppose that a function $f$ satisfies on an interval
$\mathcal{I}$ the Laguerre inequality
 \begin{equation}
 \label{lagineq}
 f'^2-f f'' >0,
 \end{equation}
 and a differential equation
\begin{equation}
\label{difeqsz}
 f''-2A(x)f'+B(x)f=0,
 \end{equation}
  where $A \in \mathbb{C}, \; \; B \in \mathbb{C}^1,$ and $B$ has at most two zeros on $\mathcal{I} .$
 Define Sonin's function by
$$
S(f;x)=f^2+\frac{f'^2}{B},
$$
 then all the local maxima of $f$ in $\mathcal{I}$ are in the interval
 defined by $B(x) >0,$ and
 $$
 sgn \left( \frac{d}{dx} S(f;x) \right) =sgn (4A B-B').
 $$
\end{lemma}
\begin{proof}
We have $0 <f'^2-f f''=f'^2-2A f f'+B f^2,$ hence $B(x) >0$
provided $f' =0.$ Finally,
\begin{equation}
\label{eqprz}
 \frac{d}{dx} \left(f^2+\frac{f'^2}{B} \right) =
\frac{4A B-B'}{B^2} \, f'^2(x),
\end{equation}
and $B(x) \ne 0$ in an interval containing all the extrema of $f$
on $\mathcal{I}.$
\end{proof}

Let us make a few remarks concerning the Laguerre inequality
(\ref{lagineq}). Usually it is stated for hyperbolic polynomials,
that is real polynomials with only real zeros, and their limiting
case, so-called Polya-Laguerre class of functions. In fact, it
holds in a more general situation. Let $L(f)=f'^2-f f'',$ defining
$\mathcal{L} =\{f(x):L(f) \ge 0 \},$ we observe that $\mathcal{L}$
is closed under linear transformations $x \rightarrow a x+b.$
Moreover, since
$$
L(f g)=f^2 L(g)+g^2 L(f),
$$
$\mathcal{L}$ is closed under multiplication. Therefore
$L(x^{\alpha} )= \alpha x^{2 \alpha -2},$ yields, in particular,
the polynomial case. Much more examples may be obtain by $L \left(
e^{f} \right)=-e^{2 f} f''.$ For our purposes it is important that
(\ref{lagineq}) holds for the function $Z_k^{\alpha , \beta }(x)$
defined by (\ref{defz}), provided $-1 \le \delta_{-1} < x <
\delta_1 \le 1,$ and $ \alpha , \beta \ge 0.$

First of all we shall establish the following claim.
\begin{lemma}
\label{unimod} Let $ Z(x)=Z_k^{\alpha , \beta }(x), \; \; k \ge 1,
\; \alpha > \beta > 0. $
 Then the global maximum of the function $S(Z; x)$
  on the interval $ [ \delta_{-1}, \, \delta_1 ]$ is
 attained at a point $x_0 =-q s - \theta \, \sqrt{(1- q^2 )(1- s^2)} \, ,$
 where $0 < \theta < \frac{2}{3} \, .$
\end{lemma}
Since the maximum of $S(Z;x)$ is also a local maximum of $Z$ by
(\ref{eqprz}), as an immediate corollary we obtain
\begin{theorem}
\label{thunimod} The global maximum of the function $Z_k^{\alpha ,
\beta }(x), \; \; k \ge 1, \; \alpha > \beta > 0,$ is attained at
a point $x_0$ defined in Lemma \ref{unimod}.
\end{theorem}

To prove Lemma \ref{unimod} we will need the following explicit
expressions.
 It is easy to check that $Z$ satisfies the following ODE
\begin{equation}
Z''-2A Z'+B z=0,
\end{equation}
where
$$
A(x)= - \, \frac{x^3+3 q s x^2+(2 q^2+ 2 s^2-1)x+ q s
}{2(1-x^2)d(x)} \, ,
$$
$$
B(x)=  \frac{d(x) r^2}{4(1-x^2)^2} \, + \frac{E(x)}{4(1-x^2)\,
d^2(x)},
$$
$$
E(x)=2 \, q s x^3-(1-4 q^2-4  s^2 + q^2 s^2 )x^2+ 6  q  s  x +1-
q^4- s^4 + 3 q^2 s^2 .
$$
Calculations yield
$$
D(x)=2(1-x^2)^2 d^3(x) (4 AB- B')=
$$
$$
q s x^6+(4  q^2+4 s^2-5 q^2 s^2 -1) x^5 +  q s(12- q^2- s^2+ q^2
 s^2 )x^4+
$$
$$2(1+  q^2+  s^2- 5 q^4 - 5 s^4 -5  q^2  s^2 + q^4  s^2+ q^2  s^4)x^3 -
$$
$$
- q  s (7+10 q^2+10  s^2-4  q^2 s^2 + q^4+ s^4)x^2 -
$$
$$(1+ 6
q^2 +6  s^2 - 6  q^4 -6  s^4-q^6-s^6+9  q^2 s^2+3  q^2 s^4+3  q^4
 s^2)x -
$$
$$
-3  q  s (2- q^2- s^2- q^4- s^4 +3  q^2 s^2).
$$
It is quite surprisingly that this expression does not contain
$r.$

 Lemma \ref{unimod} is an immediate corollary
of Lemma \ref{mainl} and the following claims
\begin{lemma}
\label{bbb} $B(x) \ne 0$  for $x \in ( \delta_{-1}, \, \delta_1
),$ provided $ k \ge 0, \; \alpha > \beta > 0.$
\end{lemma}
\begin{lemma}
\label{ddd} $D(x)  $  has the only zero $x_0$ in the interval $[
\delta_{-1}, \delta_1],$ provided $ k \ge 0, \; \alpha > \beta
>0.$
\end{lemma}
\begin{lemma}
\label{ddd1}
$$x_0 \in \left(-q s-\frac{2}{3} \, \sqrt{(1-q^2)(1-s^2) } \, , -q s
\right)\subset ( \delta_{-1}, \, \delta_1 ),$$ provided $ k \ge 0,
\; \alpha > \beta > 0.$
\end{lemma}
\begin{proof}[Proof of Lemma \ref{unimod}]
By Lemmas \ref{mainl} and \ref{bbb},
$$
sgn \left(\frac{d}{dx} S(Z;x) \right)=sgn (4AB-B')= sgn D(x).
$$
 We find
\begin{equation}
\label{dj} D(\delta_{\pm j})=-j \cos^3{ \tau} \cos^3 \omega  \,
\sin^4 (\tau+j \, \omega ) , \; \; \; j= \pm 1,
\end{equation}
hence $D(\delta_{-1} )>0, \;\; D(\delta_1 )<0.$  Therefore $x_0$
is the only maximum of $S(Z;x)$ on $[\delta_{-1}, \delta_1 ].$
\end{proof}
Let us prove Lemmas \ref{bbb}-\ref{ddd1}.
\begin{proof}[Proof of Lemma \ref{bbb}]
It is enough to show that $E(x)>0$ for $\delta_{-1} < x <
\delta_1.$ Mapping the interval $[\delta_{-1},\delta_1)$ onto $[0,
\infty )$ and simplifying we have
$$
\frac{2(1+x)^3}{ \cos^2 \tau \cos^2 \omega }\; E \left(
\frac{\delta_{-1}+ \delta_1 x}{1+x}\right)=10 \sin^2 (\tau+ \omega
)\; x^3+$$
$$
 \left( 5+2 \cos 2\tau +2 \cos 2 \omega - \cos2( \tau+
\omega ) \right)  x^2+$$
$$ \left( 5+2 \cos 2\tau +2 \cos 2 \omega
- \cos2( \tau- \omega ) \right) x+
  10 \sin^2 (\tau- \omega )>0.$$
This completes the proof.
\end{proof}

Let us remind that the discriminant of a polynomial
$p=\sum_{i=0}^n a_i x^i$ with the zeros $x_1,...,x_n ,$ is defined
by
$$
Dis_x p=a_n^{2n-2} \prod_{i <j}(x_i-x_j)^2 ,
$$
and can by calculated by the formula
$$
Dis_x p=(-1)^{n(n-1)/2}a_n ^{-1} Result_x (p,p'),
$$
where $Result_x (p,p')$ states for the resultant of $p$ and $p'$
in $x.$

 It will be a convenient tool to establish
positivity of the involved multivariable polynomials.  We used
Mathematica to find the required resultants.
\begin{proof}[Proof of Lemma \ref{ddd}]
 First we map the interval
$[\delta_{-1},\delta_1)$ onto $[0, \infty )$ by considering
$$
\frac{(1+x)^6}{  (1-q^2)^{3/2}(1-s^2)^{3/2}} \, D\left(
\frac{\delta_{-1}+\delta_1 x}{1+x}\right)=\sum_{i=0}^6 \left(v_i+
\sqrt{(1- q^2)(1- s^2} \, u_i \; \right)x^i .
$$
We will show that the sign pattern of this polynomial is
$(+++0---).$ Hence by Descartes' rule of signs it has just one
positive zero. Since by (\ref{dj})
$$
D( \delta_j )\ne 0 , \; \; \; j= \pm 1,
$$
this implies the required claim.\\
We have the following explicit expressions
$$
v_0=15 \left((s^2-q^2)^2 +8 q^2 s^2(1-q^2)(1-s^2) \right) > 0,
$$
$$
u_0=-60q s (q^2+s^2-2q^2 s^2 ) < 0;
$$
$$
v_1=12 \left( 4(1-q^2)(1-s^2)(q^2+s^2+2q^2
s^2)+3(s^2-q^2)^2\right) > 0,
$$
$$
u_1=- 96 q s (1-q^2 s^2) < 0;
$$
$$
v_2= 27(s^2-q^2)^2+40(1-q^2)(1-s^2)(q^2+s^2)+8(2-q^2-s^2+q^2 s^2)
> 0,
$$
$$
u_2=- 4q s(16-7q^2-7s^2-2q^2 s^2) < 0;
$$
$$
v_3=u_3=0;
$$
$$
v_4=-v_2 , \; \; \;  u_4= u_2;
$$
$$
v_5=-v_1, \; \; \; u_5=u_1;
$$
$$
v_6=-v_0, \; \; \; u_6=u_0;
$$
Thus, to prove the claim it is left to show that
$$w_i=v_i^2-(1-q^2)(1-s^2)u_i^2 > 0, \; \; i=0,1,2.$$
We find
$$
w_0=225 (s^2-q^2)^4 > 0;
$$
$$
w_1=144 (s^2-q^2)^2 \left(8(1-q^2 s^2)(2-q^2-s^2)
+(s^2-q^2)^2\right) > 0;
$$
$$
w_2=  \left(729-864 {\overline s}+160{\overline s}^2
\right){\overline q}^4+4 \left(135-104{\overline s}+16 {\overline
s}^2 \right) {\overline s} \, {\overline q}^3+
$$
$$
2\left(11-208{\overline s}+80 {\overline s}^2 \right) {\overline
s}^2 \, {\overline q}^2 +108 \left(5-8 {\overline s}
\right){\overline s}^3 \, {\overline q}+729{\overline s}^4 ,
$$
where
\begin{equation}
\label{zamen} {\overline q}= \sqrt{1-q} \,  ,  \; \; \; {\overline
s}= \sqrt{1-s} \, , \; \; \; 0 <{\overline s}<{\overline q} \le 1
.
\end{equation}
To demonstrate that $w_2 >0$ we calculate
$$2^{-32} \cdot 3^{-8} Dis_{\overline q} \; w_2 =$$
$$
{\overline s}^{12}(1-{\overline s})^2(45-10 {\overline
s}+{\overline s}^2)^2(360-200{\overline s}-84{\overline
s}^2-24{\overline s}^3-25 {\overline s}^4) \ne 0,$$ for
${\overline s } \in (0,1).$ Therefore $w_2$ has the same number of
real zeros for any ${\overline s} \in (0,1).$ Since ${\overline s}
<{\overline q},$ and for sufficiently small ${\overline s}>0$ we
have $w_2=729{\overline q}^4+O({ \overline s} )
>0,$ then $w_2$ has no real zeros in the region $0 <{\overline s}<{\overline
q} < 1.$ This completes the proof.
\end{proof}

\begin{proof}[Proof of Lemma \ref{ddd1}]
By Lemma \ref{ddd} it is enough to show
\begin{equation}
\label{minus} D (-q s)< 0.
\end{equation}
\begin{equation}
\label{plus} D \left( -q s-\frac{2}{3} \,  \sqrt{(1-q^2)(1-s^2)}
\, \right)  > 0.
\end{equation}
We have
$$D(-q s)= - q s (1-q^2)^2 (1-s^2)^2 \left(5+q^2+s^2 - 7 q^2 s^2 \right) <0. $$
proving (\ref{minus}). \\
To prove (\ref{plus}) we
 use again the change of variables (\ref{zamen})
obtaining
\begin{equation}
\label{eqp1p2} \frac{729}{5{\overline q}^{3/2}{\overline s}^{3/2}}
\, D \left( -q s-\frac{2}{3} \, \sqrt{(1-q^2)(1-s^2)} \, \right)=6
p_1- \sqrt{{\overline q} \, {\overline s} (1-{\overline q})(1-
{\overline s})} \, p_2 ,
\end{equation}
where
$$
p_1 =1223{\overline q} \, {\overline s}(1-{\overline
q})(1-{\overline s})+ 189({\overline q}-{\overline
s})^2+{\overline q} \, {\overline s} (93-88 {\overline q} \,
{\overline s}) >0,
$$
$$
p_2=3942 {\overline q}+3942 {\overline s}-6815{\overline q}\,
{\overline s} >0.
$$
Multiplying (\ref{eqp1p2}) by the conjugate yields
$$
h=36p_1^2-{\overline q}{\overline s}(1-{\overline q})(1-{\overline
s})p_2^2= (9-5{\overline s})(142884-43200{\overline
s}-21500{\overline s}^2+13625{\overline s}^3){\overline q}^4-$$
$$5 (555012-221688{\overline s}+127205{\overline
s}^2-46025{\overline s}^3){\overline s}\, {\overline q}^3+36
(87988+30790{\overline s}+625{\overline s}^2){\overline
s}^2{\overline q}^2-
$$
$$
4860 (571+227 {\overline s}){\overline s}^3 {\overline q}+1285956
{\overline s}^4.
$$
and
$$
2^{-4} \cdot 3^{-15} \cdot 5^{-3} Dis_{\overline s} \; h=$$
$$
{\overline s}^{12}(1-{\overline s})^2 (64459584 - 111438880
{\overline s} + 47706875 {\overline s}^2)^2 (385494997824 +
449720822304 {\overline s} -$$
$$ 674713759120 {\overline s}^2 +
      240459844600 {\overline s}^3 -
      18815866125 {\overline s}^4 +
       282081250 {\overline s}^5 -
      79028125 {\overline s}^6) \ne 0,
$$
for ${\overline s } \in (0,1).$ \\
Since $h>0$ for $s \rightarrow 0^{(+)}$ we conclude that $h >0$,
thus proving (\ref{plus}).
\end{proof}

\section{Bounds in the oscillatory region}
\label{secosc}
 It will be convenient to introduce the parameter
$\rho=2k+ \alpha+ \beta =r-1$ and two functions
$$
\mu_j(x)=\frac{\sqrt{(\rho^2-\eta^2)(\rho^2-\sigma^2)}+j(x
\rho^2+\eta \sigma) }{\rho} \, , \; \; \; j= \pm 1.
$$
The Christoffel function in the standard normalization maybe
written as
$$\frac{2k!\, \Gamma (k+ \alpha+ \beta+1)}{(2k+ \alpha+ \beta)\Gamma (k+ \alpha )
\Gamma (k+  \beta )}  \left(  P_{k-1}^{(\alpha , \beta )} (x)
\frac{d}{dx} P_k^{(\alpha , \beta )} (x) - P_k^{(\alpha , \beta )}
(x) \frac{d}{dx} P_{k-1}^{( \alpha , \beta )} (x)\right) =
$$
$$
 \sum_{i=0}^{k-1} \frac{(2i+ \alpha+ \beta+1) i! \, \Gamma(i+ \alpha+
\beta +1)}{\Gamma(i+ \alpha +1)\Gamma(i+ \beta +1)} \left(
P_i^{(\alpha , \beta )} (x)\right)^2 >0.$$

Applying the identity
$$
2(k+ \alpha )(k+ \beta ) P_{k-1}^{(\alpha , \beta )} (x) =$$
$$k \left(\beta- \alpha +(2k+ \alpha+ \beta )x \right)
P_k^{(\alpha , \beta )} (x)+ (1-x^2)(2k+ \alpha+ \beta ) \,
\frac{d}{dx} P_k^{(\alpha , \beta )} (x),
$$
and (\ref{difjp})  one finds
$$
P_{k-1}^{(\alpha , \beta )} (x) \frac{d}{dx} P_k^{(\alpha , \beta
)} (x) - P_k^{(\alpha , \beta )} (x) \frac{d}{dx} P_{k-1}^{(
\alpha , \beta )} (x) =$$ $$\frac{\rho}{2(\rho^2-\eta^2)} \, W(x),
$$
where
$$
W(x)=( \rho^2-\sigma^2 )y^2-4( \eta+\sigma x)y y'+4(1-x^2)y'^2 =$$
$$
\frac{\mu_{-1}(x) \mu_1 (x)}{1-x^2} \, y^2+
\frac{\left((\eta+\sigma x)y-2(1-x^2)y' \right)^2}{1-x^2}>0.
$$
Thus, we obtain
\begin{equation}
\label{ozyw}
 y^2 < \frac{1-x^2}{\mu_{-1}(x) \mu_1 (x)} \, W(x),
\end{equation}
provided $\mu_{-1}(x) \mu_1 (x) >0.$\\
 To estimate $W(x)$ we will consider the expression
$$
W'-z(x)W
$$
which with help of (\ref{difjp}) can be written as a quadratic
$U=A y^2+B y y'+C y'^2,$ in $y$ and $y'.$ We choose $z(x)$ in such
a way that the discriminant of $U$ vanishes, namely
$$
z_j (x)=  \frac{d}{dx}\, \left( \ln \frac{ \mu_j
(x)}{(1-x)^{\alpha+1}(1+x)^{\beta+1}}  \right), \; \; \; j=\pm 1.
$$
For such a choice of $z$ the sign of $U(x)$ coincides with the
sign of
$$C =- j \, \frac{4(1-x^2) \rho}{m_j(x)} \, , \; \; \; j= \pm 1.$$
Thus, by $W >0,$ in the region
\begin{equation}
\label{defJ} \mathcal{J} = \left( \gamma_{-1} , \gamma_1 \right),
\; \; \; \gamma_j =\frac{j \,
\sqrt{\rho^2-\eta^2)(\rho^2-\sigma^2)} -\eta \sigma}{\rho^2} \, ,
\end{equation}
defined by
$$
m_j(x)>0 \, , \; \; \; j= \pm 1,
$$
we have
$$
z_{-1}(x) < \frac{W'}{W} < z_1 (x).
$$
Solving those inequalities for $ x \in \mathcal{J},$ with initial
conditions given in a point $ \xi \in \mathcal{J},$  one obtains
\begin{equation}
\label{eqozwprom}
 \frac{(1-x)^{\alpha +1} (1+x)^{\beta +1}}{(1-\xi)^{\alpha +1} (1+\xi)^{\beta +1}} \, W(x) \ge
\left\{
\begin{array}{cc}
\frac{\mu_1
(x)}{\mu_{1}(\xi)} \, W(\xi), & x \le \xi ,\\
& \\
\frac{\mu_{-1} (x)}{\mu_{-1}(\xi)} \, W(\xi), & x \ge \xi.
\end{array}
\right.
\end{equation}
We will need the value of the following integral.
\begin{lemma}
\begin{equation}
\label{integral}
 \int_{-1}^1 (1-x)^{ \alpha+1} (1+x)^{ \beta+1} W(x) d x=
 \frac{(\rho^2-\eta^2)(\rho^2-\sigma^2)}{\rho (\rho-1)} \, {\bf h}_k^2.
\end{equation}
\end{lemma}
\begin{proof}
By
$$\frac{d}{d x} P_k^{( \alpha , \beta )} (x)= \frac{k+ \alpha+
\beta+1}{2}P_{k-1}^{( \alpha+1 , \beta +1)} (x),
$$
we have
$$
\int_{-1}^1 (1-x)^{\alpha+1} (1+x)^{ \beta+1} \left(
P_{k-1}^{(\alpha , \beta )} (x) \frac{d}{dx} P_k^{(\alpha , \beta
)} (x) - P_k^{(\alpha , \beta )} (x) \frac{d}{dx} P_{k-1}^{(
\alpha , \beta )} (x)\right) \, dx =
$$
$$
\frac{k+ \alpha+ \beta+1}{2} \int_{-1}^1 (1-x)^{\alpha+1} (1+x)^{
\beta+1}P_{k-1}^{(\alpha , \beta )} (x)P_{k-1}^{(\alpha+1 ,
\beta+1 )} (x)\,  dx \, -
$$
$$
\frac{k+ \alpha+ \beta}{2} \int_{-1}^1 (1-x)^{\alpha+1} (1+x)^{
\beta+1}P_{k}^{(\alpha , \beta )} (x)P_{k-2}^{(\alpha +1, \beta+1
)} (x) \, dx .
$$
 Now the result follows from the orthogonality relation by the
repeating application of the identities
$$
(2k+ \alpha+ \beta )P_k^{(\alpha -1, \beta)}(x)=(k+ \alpha+
\beta)P_k^{(\alpha , \beta)}(x) -(k+ \beta )P_{k-1}^{(\alpha ,
\beta)}(x)
$$
$$
(2k+ \alpha+ \beta )P_k^{(\alpha , \beta-1)}(x)=(k+ \alpha+
\beta)P_k^{(\alpha , \beta)}(x) +(k+ \alpha )P_{k-1}^{(\alpha ,
\beta)}(x),
$$
in order to express $P_i^{(\alpha , \beta)}(x) $ as a sum of
Jacobi polynomials with the parameters $\alpha +1, \beta+1 .$
\end{proof}
\begin{lemma}
\label{lelemw}
 For $x \in \mathcal{J},$
\begin{equation}
\label{lemw} (1-x)^{\alpha+1} (1+x)^{\beta+1} W(x) \le \frac{\rho
\,  \sqrt{(\rho^2-\eta^2)(\rho^2-\sigma^2)}}{ \rho-1} \, {\bf
h}_k^2 .
\end{equation}
\end{lemma}
\begin{proof}
By (\ref{eqozwprom}) and (\ref{integral}) we have
$$
 \frac{(\rho^2-\eta^2)(\rho^2-\sigma^2)}{\rho (\rho-1)} \, {\bf h}_k^2 =
 \int_{-1}^1 (1-x)^{ \alpha+1} (1+x)^{ \beta+1} W(x) d x \ge
$$
$$
(1-\xi)^{ \alpha+1} (1+\xi)^{ \beta+1} W(\xi) \left(
\int_{\gamma_{-1}}^{\xi} \frac{\mu_{1} (x)}{\mu_{1}(\xi)} \, dx
+\int_{\xi}^{\gamma_1} \frac{\mu_{-1} (x)}{\mu_{-1}(\xi)} \, dx
\right)=
$$
$$
(1-\xi)^{ \alpha+1} (1+\xi)^{ \beta+1} W(\xi) \,
\frac{\sqrt{(\rho^2-\eta^2)(\rho^2-\sigma^2)}}{\rho^2} \, ,
$$
and the result follows.
\end{proof}
\begin{lemma}
\label{lepointwise} For $x \in \mathcal{J},$
\begin{equation}
\label{oscbound} (1-x)^{ \alpha} (1+x)^{ \beta}\left({\bf
P}_k^{(\alpha , \beta ) }(x)\right)^2 < \frac{
\sqrt{(1-q^2)(1-s^2)}}{ 1-q^2-s^2-2q s x-x^2} \, .
\end{equation}
\end{lemma}
\begin{proof}
Let us remaind that $\rho=r-1, \; \eta= q r, \; \sigma =s r.$
Combining (\ref{lemw}) with (\ref{ozyw}) and using
$$
\mu_{-1}(x) \mu_1 (x)=(1-x^2) \rho^2 -2\eta \sigma
x-\eta^2-\sigma^2,
$$
we obtain the following pointwise bound.
$$
 (1-x)^{ \alpha} (1+x)^{ \beta}\left({\bf
P}_k^{(\alpha , \beta ) }(x)\right)^2 \le $$
 $$\frac{\rho}{\rho-1}
\, \frac{ \sqrt{(\rho^2-\eta^2)(\rho^2-\sigma^2)}}{ (1-x^2) \rho^2
-2\eta \sigma x-\eta^2-\sigma^2 }<
 \frac{\rho}{\rho-1} \, \frac{
\sqrt{(\rho^2-\eta^2)(\rho^2-\sigma^2)}}{(1-q^2-s^2-2q s x-x^2)
r^2} \, .
$$
It is left to check that
$$ \sqrt{(1-q^2)(1-s^2)}
-\frac{\rho}{\rho-1} \, \frac{
\sqrt{(\rho^2-\eta^2)(\rho^2-\sigma^2)}}{r^2} >0.
$$
 Multiplying this by the conjugate  and
writing it down in the variables $\alpha , \; \beta , \; k'=k+1,$
one obtains an expression with nonnegative terms only.
\end{proof}
\begin{remark}
In \cite{erdel} another pointwise estimate of the order $O \left(
\frac{\sqrt{1-x^2}}{1-q^2-s^2-2q s x-x^2} \right) $ was given. The
advantage of (\ref{oscbound}) is that it is stronger for $s =1-o(
1),$ i.e. when $ \alpha \gg k.$
\end{remark}
\begin{remark} Conjecture \ref{conjik2} would imply that
(\ref{lemw}) and, consequently, (\ref{oscbound}) is sharp up to a
multiplicative constant factor. In turn, this would imply that the
bounds of Theorems \ref{thdelta} and \ref{themn} are sharp. It is
also known that similar results hold for Laguerre polynomials
\cite{krasl2}.
\end{remark}
\section{Proof of Theorems \ref{thdelta} and \ref{themn}}
\label{secproof}
 Theorem \ref{thdelta} is an easy corollary of
Theorem \ref{thunimod} and (\ref{lepointwise}). First we need the
following
\begin{lemma}
\label{lempty}
 Suppose that $k \ge 1, \; \; \alpha > \beta >
0.$ Then
$$[-q s- \frac{2}{3} \, \sqrt{(1-q^2)(1-s^2)} \, , -q s] \subset \mathcal{J}. $$
\end{lemma}
\begin{proof}
To prove the claim it is enough to check that
$$
p(x)= (x- \gamma_{-1})(\gamma_1-x) >0,
$$
for $x=-q s$ and $x=-q s- \frac{2}{3} \, \sqrt{(1-q^2)(1-s^2)}\,
.$\\
 Straightforward calculations yield that $\rho^2 r^4 p(-q s)$
written in variables
$k'=k+1, \alpha , \beta$ is a polynomial without negative terms.\\
Similarly,
$$\rho^2 p(-q s- \frac{2}{3} \,
\sqrt{(1-q^2)(1-s^2)} \;)= r^{-4} p_1(k',\alpha, \beta)+
\frac{4}{3} \, q s (2r-1) \sqrt{(1-q^2)(1-s^2)} \, ,$$ where
$p_1(k',\alpha, \beta)$ is a polynomial in $k'=k+1, \alpha, \beta$
without negative terms.
\end{proof}

\begin{proof}[Proof of Theorem \ref{thdelta} ]
Let $\mathcal{I}=[-q s- \frac{2}{3} \, , \sqrt{(1-q^2)(1-s^2)} \,
-q s].$ Since
$$
(x-\delta_{-1})(\delta_1-x)=1-q^2-s^2-2q s x-x^2,
$$
 by Lemma \ref{lempty}, Theorem \ref{thunimod} and
(\ref{oscbound}), we get
$$
\max_{x \in [\delta_{-1},\delta_1]} Z^2(x)= \max_{x \in
\mathcal{I}} Z^2(x) < \max_{x \in \mathcal{I}}
\sqrt{\frac{(1-q^2)(1-s^2)}{ 1-q^2-s^2-2q s x-x^2}} = \frac{3
\sqrt{5}}{5} \,.
$$
This completes the proof of Theorem \ref{thdelta}.
\end{proof}

To deduce Theorem  \ref{themn} from Theorem \ref{thdelta} we will
need the following result which has been established in
\cite{krasjac}.
\begin{theorem}
\label{grmax} Suppose that $k \ge 6, \; \; \alpha \ge \beta \ge
\frac{1+ \sqrt{2} }{4} \, .$ Let $x$ be a point of a local
extremum of
$$
 (1-x)^{\alpha+1/2} (1+x)^{\beta+1/2} \left(P_k^{ \alpha, \beta} (x)\right)^2.$$
 Then $x \in \left(N'_{-1}, N'_1
 \right),$
 where
 \begin{equation}
 \label{eqN}
 N'_j= j \left(\cos (\tau'+j \omega)-
 \frac{3}{10} \, \left( \frac{\sin^4 (\tau'+j \omega)}{2 \cos \tau' \cos \omega} \right)^{1/3}
 r^{-2/3} \right),
 \end{equation}
$$
\sin \tau'= \frac{\alpha+\beta+1}{2k+ \alpha+ \beta+1} \, , \; \;
0 < \tau'  < \frac{\pi}{2} \, .
$$
 \end{theorem}
We have to restate (\ref{eqN}) in terms of $\tau.$
\begin{lemma}
\label{lvklint}
\begin{equation}
\label{vklint} (N'_{-1}, N'_1) \subset(N_{-1},N_1) \subset
(\delta_{-1}, \delta_1),
\end{equation}
where
\begin{equation}
 \label{eqN'}
 N_j= j \left(\cos (\tau+j \omega)- \frac{5}{17}  \,
  \left( \frac{\sin^4 (\tau+j \omega)}{2 \cos \tau \cos \omega} \right)^{1/3}
 r^{-2/3} \right),
 \end{equation}
 provided $k \ge 6.$
 \end{lemma}
\begin{proof}
Since $\sin \tau'=  \sin \tau +\frac{1}{r} \,,  \; \; 0<\tau ,
\tau' < \frac{\pi}{2} \,  ,$ then $\tau'
> \tau$ and $0 <\tau' \pm \omega < \pi.$ Hence
$$
\cos (\tau'+ j \omega )< \cos (\tau+ j \omega ), \; \; j=\pm 1,
$$
and
$$
[N_{-1} , N_1] \subset [- \cos (\tau'-\omega ), \cos (\tau'+\omega
)] \subset [- \cos (\tau-\omega ), \cos (\tau+\omega )]= [
\delta_{-1} , \delta_1 ].
$$
We also have for $k \ge 6,$
$$
\frac{\cos^2 \tau'}{\cos^2 \tau}=1- \frac{2 \alpha+2
\beta+1}{(2k+1)(2k+2 \alpha+2 \beta+1)} > 1- \frac{1}{2k+1} \ge
\frac{12}{13} \, .
$$
Thus
$$
\frac{\sin^4 (\tau'+ j \omega )}{\cos \tau'} >
\sqrt{\frac{12}{13}} \, \frac{\sin^4 (\tau+ j \omega )}{\cos
\tau},
$$
and as $\frac{3}{10}  \cdot
 \left(\frac{12}{13}
 \right)^{1/6}>\frac{5}{17}$ this implies $(N'_{-1}, N'_1) \subset(N_{-1},N_1).$
\end{proof}
\begin{proof}[Proof of Theorem \ref{themn} ]
To prove Theorem \ref{themn} we will bound $ M(x)=M_k^{\alpha ,
\beta }(x) $ by $Z(x).$

 Set $\epsilon_j = \left( \frac{\sin^4 (\tau+j
\omega)}{2 \cos \tau \cos \omega} \right)^{1/3}
 r^{-2/3} .$
 First we notice that
 $$
 \frac{\epsilon_j^3}{\cos\tau^3 \cos^3 \omega}=$$
 $$ \frac{(\tan \tau+ \tan \omega )^4}{2 r^2}
  <\frac{8 \tan^4 \tau }{r^2}= \frac{8(\alpha+ \beta)^2}{(2k+1)^2 (2k+2 \alpha+2 \beta+1)^2}
 <   \frac{2}{(2k+1)^2} \le \frac{2 }{169} \, .
 $$
Hence
$$
2 \cos \tau \cos \omega- \frac{5 \epsilon_j }{17}
>\left(2-\frac{5}{17}
  \left( \frac{2}{169}
\right)^{1/3} \right)\cos \tau \cos \omega > \frac{27}{14} \cos
\tau \cos \omega  .
$$

 By Theorem \ref{grmax}
and Lemma \ref{lvklint} we have
$$
\max_{x \in [-1,1] } M (x)=\max_{x \in [N_{-1},N_1] } M (x) =
\max_{x \in [N_{-1},N_1] } \sqrt{\frac{1-x^2}{1-q^2-s^2-2q s
x-x^2}} \; Z^2(x) <
$$
$$
\max_{x \in [N_{-1},N_1] } \frac{3}{5}\,
\sqrt{\frac{5(1-x^2)}{1-q^2-s^2-2q s x-x^2}} <  \max_{j= \pm 1}
\frac{3}{5}\, \sqrt{ \frac{5 \sin^2 (\tau+j \omega)}{1-q^2-s^2-2q
s N_j -N_j^2}} \, =
$$
$$
\frac{3 \sqrt{17} \, \sin (\tau+j \omega )}{5 \, \sqrt{\epsilon_j
(2 \cos \tau \cos \omega-\frac{5}{17} \, \epsilon_j ) }} <
\sqrt{\frac{238}{75}} \,
  \frac{\sin(\tau+ j \omega )}{ \sqrt{\epsilon_j
\cos \tau \cos \omega }}=
$$
$$
\sqrt{\frac{238}{75}} \, \left( \frac{r\, \sin (\tau+j \omega )
}{\cos \tau \cos \omega } \right)^{1/3} =\sqrt{\frac{238}{75}} \,
r^{1/3} ( \tan \tau + \tan \omega )^{1/3} < \frac{9}{4} \, (r \tan
\tau )^{1/3} =
$$
$$
\frac{9}{4} \, \left( \frac{( \alpha+ \beta )^2 (2k+ \alpha+ \beta
+1)^2 }{(2k+1)(2k+ 2 \alpha+2 \beta+1)}\right)^{1/6} \le
\frac{9}{4} \, \left( \frac{4 \alpha^2  \,(2k+ 2 \alpha +1)^2
}{(2k+1)(2k+ 4 \alpha+1)}\right)^{1/6}
 < $$
 $$3  \alpha^{1/3} \left(1+
\frac{\alpha}{k} \right)^{1/6}.
$$
This completes the proof.
\end{proof}


\end{document}